\documentstyle[12pt]{article}
\catcode`\@=11
\@addtoreset{equation}{section}
\def\theequation{\thesection.\arabic{equation}}
\catcode`\@=12
\newtheorem{Theorem}{Theorem}[section]
\newtheorem{Definition}{Definition}[section]
\newtheorem{Proposition}{Proposition}[section]
\newtheorem{Lemma}{Lemma}[section]
\newtheorem{Corollary}{Corollary}[section]
\newtheorem{Note}{Note}[section]
\title{On Existences Of Periodic Orbits For Hamilton Systems\thanks{Project 19871044 Supported by NSF}}

\author{Renyi Ma\\
Department of Mathematics \\
Tsinghua University \\
Beijing, 100084\\
People's Republic of China\\
rma@math.tsinghua.edu.cn}

\date { }

\begin{document}
\textwidth=125mm
\textheight=185mm
\parindent=8mm
\frenchspacing
\maketitle

\begin{abstract}
In this article, we prove that either there 
exists at least one periodic orbit
of Hamilton vector field on a given 
energy hypersurface in $R^{2n}$ or there 
exist at least two periodic orbits on the near-by 
energy hypersurface in $R^{2n}$.
The more general results 
are also obtained.
\end{abstract}


\noindent{\bf 2000MR Subject Classification} 32Q65,53D35,53D12

\section{Introduction and results}

Let $\Sigma$ be a smooth closed oriented manifold of dimension
$2n-1$ in $R^{2n}$, here $(R^{2n},\omega _0)(\omega _0=\sum _{i=1}^ndx_i\wedge dy_i)$ is 
the standard symplectic space. Then 
there exists a unique 
vectorfield $X_\Sigma $ so called Hamilton vector field defined
by
$i_{X_{\Sigma }}\omega _0|\Sigma \equiv 0$.
A periodic Hamilton orbit in 
$\Sigma $ is a smooth path
$x:[0,T]\to \Sigma ,T>0$
with $\dot x(t)=X_{\Sigma }(x(t)) \ for \ t\in(0,T)$ and $x(0)=x(T)$.
Seifert in \cite{se} raised the following conjecture: 

{\bf Conjecture}(see\cite{se}). Let $\omega _0$ be the standard 
symplectic form 
on $R^{2n}$. Let $\Sigma $ be 
a closed $(2n-1)-$hypersurface  
in $R^{2n}$. Then 
there is a closed Hamilton periodic orbit in $\Sigma $.

Rabinowitz \cite{ra} and Weinstein \cite{we1,we2} 
proved that if $S$ is starshape resp. convex, 
$SC$ holds. 
Weinstein conjectured $SC$ holds for 
the hypersurface of contact type in general symplectic manifold($WC$). 
In $R^{2n}$, $SC$ implies $WC$.
Viterbo \cite{vi} proved the 
$WC$ in $R^{2n}$. By the Viterbo's work, Hofer and Zehnder in \cite{hz} 
proved the near-by $SC$ holds. Struwe in \cite{st} proved that $SC$ almost holds by modifing Hofer-Zehnder's work. 
Ginzburg in \cite{gi} and Herman in \cite{he} gave a counter-example \cite{he} 
for the $SC$.
After Viterbo's work, many results 
were obtained 
by using variational method or Gromov's 
$J-$holomorphic curves via 
nonlinear Fredholm alternative, see 
\cite{fhv,ho,hv1,hz,ma1,ma2,ma3,ma4} etc.

Let $(M,\omega )$ be a symplectic manifold and 
$h(t,x)(=h_t(x))$ a compactly supported smooth function on 
$M\times [0,1]$. Assume that the segment 
$[0,1]$ is endowed with time coordinate $t$. For every function 
$h$ define the $(time-dependent)$ $Hamiltonian$ 
$vector$ $field$ $X_{h_t}$ by the equation:
\begin{equation}
dh_t(\eta )=\omega (\eta ,X_{h_t}) \ \ for \ every \ \eta \in TM
\end{equation}
The flow $g^t_h$ generated by the field 
$X_{h_t}$ is called $Hamiltonian $ $flow$ and its time one map $g_h^1$ is called $Hamiltonian $ $diffeomorphism $.

Now assume that $H$ be a time independent smooth function on 
$M$ and $X_H$ its induced vector field.

Let $(M,\omega )$ be a symplectic manifold.
Let $J$ be the almost complex structure tamed by $\omega $, i.e., 
$\omega (v, Jv)>0$ for $v\in TM$. Let 
${\cal {J}}$ the space of all tame almost complex structures. 
 
\begin{Definition}
Let 
$$s(M,\omega ,J)=\inf \{ \int _{S^2}f^*\omega >0 |f:S^2\to M \ is \ J-holomorphic \} $$
\end{Definition}

\begin{Definition}
Let 
$$s(M,\omega )=\sup _{J\in {\cal {J}}}l(M,\omega ,J)$$
\end{Definition}

Let $W$ be a Lagrangian submanifold in $M$, i.e., 
$\omega |W=0$. 

\begin{Definition}
Let 
$$l(M,W, \omega )=\inf \{ |\int _{D^2}f^*\omega |>0 |f:(D^2,\partial D^2)\to (M,W) \} $$
\end{Definition}

\begin{Theorem}
Let $(M,\omega )$ be a closed  
compact symplectic manifold or a manifold 
convex at infinity and 
$M\times C$ be a symplectic manifold with 
symplectic form $\omega \oplus \sigma $, 
here $(C,\sigma )$ standard 
symplectic plane. Let 
$2\pi r_0^2<s(M,\omega )$ and 
$B_{r_0}(0)\subset C$ the closed ball with radius 
$r_0$. 
Assume that $H$ be a time independent smooth function on 
$M\times C$ and $X_H$ its induced vector field. 
If 
$\Sigma =H^{-1}(c)$ be a smooth hypersurface in $M\times B_{r_0}(0)$, $X_H$ its Hamilton vector 
field, then 
either there exists at least one periodic orbits of $X_H$ on $\Sigma $ or 
there exists at least two periodic orbits of $X_H$ on $\Sigma '=H^{-1}(c')$, 
$c'$ is close to $c$ as one wants. 
\end{Theorem}

\begin{Theorem}
Let $(M,\omega )$ be an exact symplectic manifold, i.e., 
$\omega =d\alpha $ for some $1-form$ $\alpha $.  
Assume that $H$ be a time independent smooth function on 
$M$ and $X_H$ its induced vector field. 
If 
$\Sigma =H^{-1}(c)$ be a smooth compact hypersurface in $M$ and there exists a 
Hamiltonian diffeomorphism $h$ such that 
$h(\Sigma )\cap \Sigma =\emptyset $,  
then 
either there exists at least one periodic orbits of $X_H$ on $\Sigma $ or 
there exists at least two periodic orbits of $X_H$ on $\Sigma '=H^{-1}(c')$, 
$c'$ is close to $c$ as one wants. 
\end{Theorem}
The near-by and almost 
existence results or the finity of Hofer-Zehnder's for bounded set in the above 
symplectic manifolds can be obtained as in \cite{hzb}. For example, one has 
\begin{Corollary}
Let $M$ be any open manifold and $(T^*M,d\alpha )$ be its cotangent bundle.  
Assume that $B$ is bounded set in 
$T^*M$, then the Hofer-Zehnder capacity
$C_{HZ}(B)$ is finite. 
\end{Corollary}

Theorem1.1-1.2 was reported in 
the proceedings of 
the international conferrence on 
``Boundary Value Problems, Integral Equations, And Related Problems''(5-13 August 2002); 
``ICM2002-Beijing Satellite Conference on 
Nonlinear Functional Analysis, August 14-18,2002 Taiyuan.
Theorem 1.1-1.2
can be generalized to the products of symplectic manifolds in Theorem 1.1-1.2. The proofs 
of Theorem1.1-1.2 is close as in \cite{ma3,mo}. 
Here we flow the Monke's method in \cite{mo}.
If a $(n-1)$-dimensional submanifold $\cal L$ in 
$\Sigma$ satisfying that 
$\cal L$ is transversal to 
the hamilton vector field $X_H$ and 
$\omega _0|{\cal {L}}=0$ and 
$\omega _0|\pi _2(M,{\cal {L}})=0$, 
then we call 
${\cal {L}}$ the Hamilton-Legendre isotropic submanifold. 
A Hamilton-Arnold chord in 
$\Sigma $ is a smooth path
$x:[0,T]\to \Sigma ,T>0$
with $\dot x(t)=X_{\Sigma }(x(t)) \ for \ t\in(0,T)$ and $x(0),x(T)\in {\cal {L}}$.
Then one can also
prove the Hamilton's chord  almost existence results as Theorem1.1-1.2.

\section{Lagrangian Non-Squeezing}

\begin{Theorem}
(\cite{po})Let $(M,\omega )$ be a closed  
compact symplectic manifold or a manifold convex at infinity and 
$M\times C$ be a symplectic manifold with 
symplectic form $\omega \oplus \sigma $, 
here $(C,\sigma )$ standard 
symplectic plane. Let 
$2\pi r_0^2<s(M,\omega )$ and 
$B_{r_0}(0)\subset C$ the closed disk with radius 
$r_0$. If $W$ is a close Lagrangian manifold in $M\times B_{r_0}(0)$, then 
$$l(M,W, \omega )<2\pi r_0^2$$
\end{Theorem}
This can be considered as an Lagrangian version of Gromov's symplectic non-squeezing(\cite{gro}).

\begin{Corollary}
(Gromov\cite{gro})Let $(V',\omega ')$ be an exact symplectic manifold with 
restricted contact boundary and $\omega '=d\alpha ' $. Let 
$V'\times C$ be a symplectic manifold with 
symplectic form $\omega '\oplus \sigma =d\alpha 
=d(\alpha '\oplus \alpha _0)$, 
here $(C,\sigma )$ standard 
symplectic plane. 
If $W$ is a close exact Lagrangian submanifold, then 
$l(V'\times C,W, \omega )=\infty $, i.e., there does not exist any 
close exact Lagrangian submanifold in $V'\times C$.
\end{Corollary}

\begin{Corollary}
Let $L^n$ be a close Lagrangian in $R^{2n}$ and $L(R^{2n},L^n,\omega )=2\pi r_0^2>0$, then 
$L^n$ can not be embedded in $B_{r_0}(0)$ as a Lagrangian submanifold. 
\end{Corollary}

\section{Construction Of Lagrangian}

\subsection{First Case: no periodic orbit}

Let $(V,\omega )=(R^{2n}\times R^{2n},\omega _0\ominus \omega _0)$ 
be the standard symplectic vector space, here 
$\omega _0=d\lambda _0=d({{1}\over {2}}(x_idy_i-y_idx_i))$.
Let $\Sigma $ be a oriented closed 
hypersurface in 
$R^{2n}$. 
Let 
${\cal L}=\{ (\sigma ,\sigma )|\sigma \in \Sigma \subset 
R^{2n}\}$ 
be a closed 
isotropic submanifold 
contained in $
(\Sigma ',\omega )=(\Sigma \times \Sigma,\omega _0\ominus \omega _0)$, i.e., $Q^*\omega |_{{\cal L}}=0$.
Since $\Sigma $ is oriented in $R^{2n}$, 
the normal bundle of $\Sigma $ is trivial. So, 
the tubular neighbourhood $Q_\delta (\Sigma )$ 
of $\Sigma $ is foliated by $\Sigma $.
We now define a Hamiltonian as follows. 
Define $Q_\delta (\Sigma )
=\cup _{|t|\leq \delta }\psi _{t}(\Sigma )$, 
here $\psi _t$ is the flow of the normal vector field of $\Sigma $.
Define $h(x)=t(x)$, $t(x)$ is the arrival time 
of $\psi _t$ from $x$ to 
$\Sigma $ for $x\in Q_\delta (\Sigma )$. 
Let $H(\sigma _1,\sigma _2)=h(\sigma _2)$ on 
$\bar Q_\delta =Q_\delta (\Sigma )\times Q_\delta (\Sigma )$ 
and $X_H$ be its Hamilton vectorfield.  
Let $\eta _s$ be the Hamilton flow on $\bar Q_\delta $
induced by $X_H$. 
Let $S$ be a smooth compact hypersurface in $R^{4n}$ 
which intersects the hypersurface $E=R^{2n}\times \Sigma $ transversally 
and contains 
${\cal L}$.
Furthermore, 
$S$ is transversal to 
the hamilton vector field 
$X_H$. 
Let $s(x)$ be the time to $S$ of hamilton flow $\eta _s$ which is well defined 
in the neighbourhood $U_{\delta _1}(S)$. Let $X_S$ be the Hamilton vector field 
of $S$ on $U_{\delta _1}$ and $\xi _t$ be the flow of $X_S$ defined on $U_{\delta _1}$. 
Since $X_S(H)=\{ s,H\}=-\{ H,s\}=-X_H(s)=-1\ne 0$, here $\{ ,\}$ is the Poisson
bracket, so $X_S$ is transversal to $E$.
Then, there exists $\delta _0$ such that 
$\xi _t(x)$ exists for any $x\in U_{\delta _0}(S\cap E)$, $|t|\leq 100\delta _0$. Let 
$L_0=\cup _{t\leq \delta _0}\xi _t({\cal {L}})$. 
One has 
\begin{Lemma}
$L_0$ is a Lagrangian submanifold in $(R^{4n},\omega )$. 
\end{Lemma}
Proof.
Let
\begin{eqnarray}
&&F:[-\delta _0,\delta _0]\times {\cal {L}}\to \bar Q_\delta \cr 
&&F(t,l)=\xi _t(l)). 
\end{eqnarray}
Then,
\begin{eqnarray}
F^*\omega &=&F^*_t\omega +i_{X_S}\Omega \wedge dt\cr 
&=&0-(ds|S)\wedge dt=0
\end{eqnarray}
This checks that $L_0=F([-\delta _0,\delta _0])\times {\cal {L}})$ is Lagrangian submanifold.

\begin{Lemma}
Let $M=E\cap S$ and $\omega _M=\omega |M$. 
Then there exist $\delta _0>0$ and a neighbourhood $M_0$ of $\cal L$ in $M$ such that 
$G:M_0\times [-\delta _0,\delta _0]\times [-\delta _0,\delta _0]\to R^{4n}$ defined by  
$G(m,t,s)=\eta _s(\xi _t(m))$ is an symplectic embedding.
\begin{eqnarray}
G^*\omega =\omega _M+dt\wedge ds
\end{eqnarray}
\end{Lemma}
Proof. It is obvious.

Let $U_T=G(M_T\times [-\delta _0,\delta _0]\times [-T,T])$. 
If there does not exist periodic solution 
in $Q_\delta (\Sigma )$ and $M_T$ is a very small neighbourhood of ${\cal {L}}$, 
then  $s(x)$ and $X_S$ is well defined on 
$U_T$. 
Therefore, there exists $\delta _T$ such that 
the flow $\xi _t(x)$  of $X_S$ exists for 
any $x\in U_T$, $|t|\leq 100\delta _T$. 
Let $U_0=U_{\delta _0}(S\cap E)\cap U_T$, $U_k=\eta _{k\delta _0}(U_0)\subset U_T$, $k=1,...,k_T$.  
Let $X_k=X_S|U_k$. 
Let $\bar X_k=\eta _{k\delta _0*}X_0$.  
We claim that 
$\bar X_k=X_k$. 
Since $s(x)$ and $H(x)$ is defined 
on $U_T$ and $\{ H,s\}=1$, so $\xi _t(\eta _s(x))=\eta _s(\xi _t(x))$ for $x\in U_T$. 
Differentiate it, we get $X_S(\eta _s(x))=\eta _{s*}X_S(x)$. Take $s=k$, 
one proves 
$\bar X_k=X_k$. 
Recall that 
the flow $\xi _t(x)$  of $X_S$ exists for 
any $x\in U_{\delta _0}(S\cap E)$, $|t|\leq 100\delta _0$. 
So, the flow $\bar \xi ^k_t(x)$  of $\bar X_k$ exists for 
any $x\in U_k$, $|t|\leq 100\delta _0$. 
Therefore, the flow $\xi ^k_t(x)$  of $X_k$ exists for 
any $x\in U_k$, $|t|\leq 100\delta _0$. This Proves that 
the flow $\xi _t(x)$  of $X_S$ exists for 
any $x\in U_T$, $|t|\leq 100\delta _0$.

\begin{Theorem}
(Long Darboux theorem)
Let $M=E\cap S$, $\omega _M=\omega |M$. Let 
$M_0$ as in Lemma 3.2. Let 
$(U'_T,\omega ')=(M_T\times [-\delta _0,\delta _0]\times [-T,T],\omega _M+dH'\wedge ds')$. 
If there does not exist periodic solution 
in $Q_\delta (\Sigma )$, then there exists a symplectic embedding 
$G:U'_T \to \bar Q_\delta $ defined by 
$G(m,H',s')=\eta _{s'}(\xi _{H'}(m))$ such that 
\begin{eqnarray}
G^*\omega =\omega _M+dH'\wedge ds'.
\end{eqnarray}
\end{Theorem}
Proof. We follow the Arnold's proof 
on Darboux's theorem in \cite{arb}. 
Take a Darboux chart $U$ on $M$, we assume 
that  
$\omega _M|=\sum _{i=1}^{2n-1}dp'_i\wedge dq'_i$. Now computing the 
Poisson brackets $\{,\}^*$ of $(p'_1,q'_1;...,p'_{2n-1},q'_{2n-1};H',s')$ for 
$G^*\omega $ on $U'_T$. Let 
$P'_i(G(p'_1,q'_1;...,p'_{2n-1},q'_{2n-1};H',s'))=P'_i(\xi _{H'}\eta _{s'}(p',q'))=p'_i$, and
$Q'_i(G(p'_1,q'_1;...,p'_{2n-1},q'_{2n-1};H',s'))=Q'_i(\xi _{H'}\eta _{s'}(p',q'))=q'_i$,
$i=1,...,2n-1$.
$\{H',s'\}^*=G^*\omega ({{\partial }\over {\partial H'}},{{\partial }\over {\partial s'}})
=\omega (G_*{{\partial }\over {\partial H'}},G_*{{\partial }\over {\partial s'}})
=\omega (X_H,X_S)=1$, 
$\{H',p_i'\}^*=G^*\omega ({{\partial }\over {\partial H'}},{{\partial }\over {\partial p_i'}})
=\omega (G_*{{\partial }\over {\partial H'}},G_*{{\partial }\over {\partial p_i'}})
=\omega (X_H,X_{P_i'})=-X_{P_i'}(H)=0$. Similary, 
$\{H',q_i'\}^*=0$. Similarly, 
$\{s',H'\}^*=1$, $\{s',p_i'\}^*=0$,$\{s',q_i'\}^*=0$. 
Note that 
$\omega =\xi _t^*\eta _s^*\omega $, so 
$\{p_i',q_j'\}^*=
G^*\omega ({{\partial }\over {\partial p_i'}},{{\partial }\over {\partial q_j'}})
=\omega (G_*{{\partial }\over {\partial p'_i}},G_*{{\partial }\over {\partial q'_j}})
=\omega (\xi _{H'*}\eta _{s'*}{{\partial }\over {\partial p'_i}},
\xi _{H'*}\eta _{s'*}{{\partial }\over {\partial q'_j}})
=
\xi ^*_{H'}\eta ^*_{s'}
\omega ({{\partial }\over {\partial p'_i}},
{{\partial }\over {\partial q'_j}})
=\omega ({{\partial }\over {\partial p'_i}},
{{\partial }\over {\partial q'_j}})=\omega _M({{\partial }\over {\partial p'_i}},
{{\partial }\over {\partial q'_j}})=\delta _{ij}$. This shows that 
the 
Poisson brackets $\{,\}^*$ is same as 
the Poisson brackets $\{,\}'$ for $\omega '$. So, 
$\omega '=G^*\omega $. This finishes the proof.

Take a disk $M_0$ enclosed by the circle $E_0$ which is 
parametrized by $t\in [0,\delta _0]$ in $(s',H')-plane$ 
such that 
$M_0\subset [-2s_0,2s_0]\times [0,\varepsilon]$ and $area (M_0)\geq 2s_0\varepsilon $. 
Now one checks that 
$L=G({\cal {L}}\times E_0)$ satisfy  
\begin{eqnarray}
\omega |L=G^*\omega ({\cal {L}}\times E_0)&=&dH'\wedge ds'|E_0=0. 
\end{eqnarray}
So, $L$ is a Lagrangian submanifold.

\begin{Lemma}
If there does not exist any periodic orbit 
in $(Q_\delta ,X_H)$, 
then $L$ 
is a close Lagrangian submanifold. Moreover 
\begin{equation}
l(V,L,\omega )=area(M_0)
\end{equation}
\end{Lemma}
Proof. 
It is obvious that  $F$ is a Lagrangian embedding.
If the circle $C$ homotopic to $C_1\subset {\cal {L}}\times s_0$ then  we compute
\begin{eqnarray}
\int _CF^*(p_idq_i)=\int _{C_1}F^*(p_idq_i)=0. 
\end{eqnarray}
since $\lambda |C_1=0$ due to $C_1\subset {\cal {L}}$ and 
$\cal L$ is ``Legendre submanifold". 

If the circle $C$ homotopic to $C_1\subset l_0\times S^1$ then  we compute
\begin{eqnarray}
\int _CF^*(p_idq_i )=\int _{C_1}F^*(p_idq_i )=n(area(M_0)). 
\end{eqnarray}
This proves the Lemma.

\subsection{Second case: single curve of periodic orbits}

\begin{Theorem}
(Long Darboux cover theorem)
Let $M_0$ be as in Lemma 3.2. 
Let 
$\bar U'_T=M_0\times [-\delta _0,\delta _0]\times [-T,T]$, 
$\omega '=\omega _M+dH'\wedge ds'$. 
Then there exists a symplectic immersion 
$G:U'_T\to \bar Q_\delta $ satisfies 
$G(m;H',s')=\xi _{H'}\eta _{s'}(m)$;
and
\begin{eqnarray}
G^*\omega =\omega _M+dH'\wedge ds'.
\end{eqnarray}
\end{Theorem}
Proof. By the proof of Theorem 3.1.

Now Assume that $H$ is ``single'', i.e., there exist only one  family of periodic orbits, more precicely, 
the periodic orbits consist of $x(t,(\sigma _0,c))$ for some $0<c<\delta $ in $Q_\delta $ such that 
$H(x(t,(\sigma ,c))=c\}$ with period $T_c$. 
Now we assume that there does not exist any peiodic orbit on 
$H^{-1}(0)$. 
Let 
$Z'_T={\cal {L}}\times [-\delta _0,\delta _0]\times [-T,T]$. 
$Z_T=G(Z_T')$. 
Let 
$\gamma (l)=L_0\cap (Z_T\setminus L_0)$ and 
$\{\gamma '_i \}_{i=1}^m=G^{-1}(\gamma )\subset Z'_T$.
We claim that one still can take a disk $M'_0$ enclosed by the circle $E'_0$ which is 
parametrized by $t\in [0,t_0]$ in $(s',H')-plane$ 
such that 
$G|E'_0$ is an embedding and 
$M'_0\subset [-2s_0,2s_0]\times [0,\varepsilon]$ with 
$0<\varepsilon <\delta _0$ and $area (M'_0)\geq 2s_0\varepsilon $. 
In fact one can draw 
a curve $E'_1$ like rectangle without bottom under the level $H'=\varepsilon $ above 
$s'-axis$ between $\gamma _1'$ and $\gamma _2'$ on 
$(s',H')$-plane. 
Let $E_1=G(E'_1)$ and $\{E'_i\}_{i=1}^n=G^{-1}(E_1)$, 
one can draw similar graph curve $F'_2$ over $s'-axis$ below $E'_2$ between $\gamma _2'$ and $\gamma _3'$ on 
$(s',H')$-plane. We do this similarly $n$ times. Then, we connect 
$E'_1$, $F'_2$,...,etc. below $\gamma '_i$ above $s'-axis$ to get a graph curve $\Gamma $. Finally, 
we close $\Gamma $ with $s'-axis $ to get $E'_0$.  
Let $L=G({\cal {L}}\times E'_0)$. 
So, $L$ is again a Lagrangian submanifold.
The Lemma 3.2 still holds in this case.

\subsection{Gromov's figure eight construction}

First we note that the construction of section 3.1 holds for any symplectic manifold.
Now let $(M,\omega )$ be an exact symplectic manifold with 
$\omega =d\alpha $. 
Let $\Sigma =H^{-1}(0)$ be a regular and close smooth 
hypersurface in $M$ and $H$ is $T-finite$. $H$ is a time-independent Hamilton function. 
Set 
$(V',\omega ')=(M\times M,\omega \ominus \omega )$.
If there does not exist any close 
orbit for $X_H$ 
in $(\Sigma ,X_H)$, one can construct 
the Lagrangian 
submanifold 
$L$ as in section 3.1, let $W'=L$.
Let $h_t=h(t,\cdot ):M\to M$, $0\leq t\leq 1$ be a Hamiltonian isotopy 
of $M$ induced by hamilton fuction $H_t$ such that $h_1(\Sigma )\cap \Sigma =\emptyset $, 
$|H_t|\leq C_0$.
Let $\bar h_t=(id,h_t)$. Then $F'_t=\bar h_t:W'\to V'$ 
be an isotopy of Lagrangian embeddings. 
As in \cite{gro}, we can use symplectic figure eight trick invented by Gromov to 
construct a Lagrangian submanifold $W$ in $V=V'\times R^2$ through the 
Lagrange isotopy $F'$ in $V'$, i.e., we have 
\begin{Proposition}
Let $V'$, $W'$ and $F'$ as above. 
Then there exists  
a weakly exact Lagrangian embedding $F:W'\times S^1\to V'\times R^2$ 
with $W=F(W'\times S^1)$ is contained in 
$M\times M\times B_R(0)$, here $4\pi R^2=8C_0$ and 
\begin{equation}
l(V',W,\omega )=area(M'_0)=A(T).
\end{equation}
\end{Proposition}
Proof. Similar to \cite[2.3$B_3'$]{gro}.

{\bf Example.} Let $M$ be an open manifold and $(T^*M,p_idq_i)$ be 
the cotangent bundle of open manifold with the 
Liouville form $p_idq_i$. Since 
$M$ is open, there exists a function $g:M\to R$ without 
critical point. The translation by 
$tTdg$ along the fibre gives a hamilton isotopy 
of $T^*M:h^T_t(q,p)=(q,p+tTdg(q))$, so
for any given compact set $K\subset T^*M$, there exists 
$T=T_K$ such that $h^T_1(K)\cap K=\emptyset $.

\section{Proof on Theorems}

Take $T_0>0$ such that $A(T_0)\geq 100\pi r_0^2$. Assume that 
on $H^{-1}(0)$ there does not exist periodic orbit with period 
$T\leq 100T_0$ and $H$ as in section 3, 
then 
by the results in section 3, we have 
a close Lagrangian submanifold 
$W=L$ or $W=F(W'\times S^1)$ contained in 
$V=M\times C\times M\times B_{r_0}(0)$ or 
$V=M\times M\times B_{r_0}(0)$. By Lagrangian non-squeezing theorem, i.e., Theorem 2.1, 
we have 
\begin{equation}
A(T_0)\leq area(M_0)=l(V,W,\omega )\leq 2\pi r_0^2.
\end{equation}
This is a contradiction. This contradiction shows 
that 
there is a periodic orbit with period 
$T\leq 100T_0$. 
This completes the proofs of theorems.

\end{document}